\documentclass[11pt]{article}

\usepackage{graphicx} 
\usepackage{xcolor}
\usepackage{epsfig}

\title{\bf Utilization of noise for the control of a class of non-linear systems}

\author{A.-M. Stoica\footnote{
{\tt adrian.stoica@upb.ro} University POLITEHNICA of Bucharest, Faculty of Aerospace Engineering, Str. Polizu, No. 1, 011063,Romania} \and I. Yaesh\footnote{Control Department, Elbit Systems, P.O.B. 1044/77, Ramat-Hasharon, 47100, Israel}}
\date{~}
\begin{document}

\maketitle

\pagestyle{myheadings}
\markboth{Author1 Author2
}{$IAAC^3$}

\bigskip

\begin{abstract}
Utilization of noise for the control of a class of non-linear systems is presented. The application of state-multiplicative noise as a mean of control is far more limited then the use of standard deterministic gains. Nevertheless, so called Stochastic Anti Resonance (SAR) with state-multiplicative noise based control, do arise in a variety of situations such as in engineering applications, physics modelling, biology, and models of visuo-motor tasks. Linear Matrix Inequalities based conditions from recent publications are reviewed, that characterize stochastic stability of such nonlinear systems applying SAR. While those results dealt with systems that are, apriori, modelled using sector bounded nonlinearities, we demonstrate that more generals systems that can be approximated as such, can be also controlled using SAR.

\end{abstract}


\section{Introduction}

Deterministic state-feedback controllers are far more popular than stochastic control achieving Stochastic Antiresonance (SAR), where 
state-multiplicative noise is applied to stabilize non linear systems. The reverse phenomenon, Stochastic Resonance (SR), has been mentioned in connection with several topics such as periodic occurrence of ice ages in \cite{Wellens}. Other cases of SR are mentioned in \cite{Wellens} and \cite{Kalasnikov}, in relation to 
a particle in a double well, animal behavior, sensory neurons and ionic channels, optical systems and so on. 
SAR in squid giant axons has been studied in \cite{Borkowsky}, where the potential for therapeutic neurological applications was pointed out. 
Anti-resonance applications may also be found in \cite{SWG}, \cite{JYK} and \cite{Kremer}. Recently the stabilization of a class of non linear systems with sector bounded non linearities has been analyzed, where Linear Matrix Inequalities (LMI) based condition have been derived \cite{syent}. Here we shortly review those conditions and their application, and we widen their scope of applications, by applying the universal approximation theorem \cite{Cyb} for systems that are not apriori modelled with sector bounded uncertainties.



\section{Review of Stability Results for Sector Bounded Nonlinearities}

Consider the following system:
\begin{eqnarray}
	\label{e1bis}
	\begin{array}{lcl}
		dx(t)& =& Ax(t)dt + F f(y(t))dt + Dx(t)d\beta (t)   \label{sys} \\
		y\left( t \right) &=& Cx\left( t \right)\\
		x(0)&=&x_0 
	\end{array}
\end{eqnarray}
where $x \in {\mathcal R^n}$ denotes the state vector, $y \in {\mathcal
	R^n}$ is the measured system output and where $\beta \left( t\right)
\in {\mathcal R} $ is a standard Wiener
process with $E\{\beta^2(t)\}=dt$ on the given probability space which is also
independent of $x_0$. The elements of $y$ are $y_i=C_ix\in {\mathcal{R}}, \, i=1,...,n$, where $C_i$ is the $i$'th row vector of $C$, namely $y_i=\sum_{j=1}^nC_{ij}x_j$ and 
the components $f_i(y_i) $ of $f(y)$ satisfy the sector conditions 
$0 \leq y_i f_i(y_i) \leq s_i y_i^2$  (\cite{Banjerdpongchai}, \cite{Lure}, \cite{Boyd}) which are equivalent to
\begin{equation}  
	\label{sector1}
	f_i(y_i) (f_i(y_i) - s_i y_i) \leq 0,\, i=1,...,n.
\end{equation}
Let us define, for the sequel $S = diag\{s_1,s_2,...,s_n\}$.  

We also assume that the derivatives of the nonlinearities are bounded, i.e. $\delta_i>0$, $i=1,\dots,n$ such that $\frac{ df_i(y_i)}{dy_i}<\delta_i$, and 
$C^TC=I$. It was shown in  \cite{syent} that if there exist $\nu\in(0,1)$, $\Lambda=diag(\lambda_1,\dots,\lambda_n)$, $\lambda_i\geq 0$, $i=1,\dots, n$ and  $\mathcal{T}=diag(\tau_1,\cdots,\tau_n),\, \tau_i\geq 0$, $i=1,\dots, n$ such that
	\begin{eqnarray}
		\label{e18}
		\mathcal{N} : = \left[\begin{array}{cc}{\mathcal{N}}_{11}(\nu,\Lambda)&{\mathcal{N}}_{12}(\nu, \Lambda,\mathcal{T})\\{\mathcal{N}}_{12}^T(\nu, \Lambda,\mathcal{T})&{\mathcal{N}}_{22}(\Lambda, \mathcal{T})\end{array}\right]<0
	\end{eqnarray}
	where $\Delta:=diag\left(\delta_1,\dots,\delta_n\right)$ and where
	\begin{eqnarray}
		\label{e19}
		\begin{array}{lcl}
			{\mathcal{N}}_{11}(\nu,\Lambda)&:=&\nu\left[A^T+A-\sigma^2(1-\nu)I\right]
			+\sigma^2C^T\Lambda \Delta C\\
			{\mathcal{N}}_{12}(\nu, \Lambda,\mathcal{T})&:=&\nu F+\left[A-\frac{\sigma^2}{2}\left(1-\frac{\nu}{2}\right)I\right]^TC^T\Lambda
			+SC^T\mathcal{T}\\
			{\mathcal{N}}_{22}(\Lambda, \mathcal{T})&:=&-2\mathcal{T}+\Lambda C F+F^TC^T\Lambda ,
		\end{array}
	\end{eqnarray}
	then the solution $x(t)\equiv 0$ of the stochastic system (\ref{e1bis}) with $D=\sigma I$ is asymptotically stable for the above system. 
	
	The derivation of the above result was applied by evaluating the infinitesimal generator (\cite{Friedman},\cite{Mao}) of the novel non-quadratic Lyapunov functional 
	
	$$
	V(x) = (x^Tx)^{\nu/2} +  \Sigma_{k=1}^n \lambda_k \int_0^{y_k}s^{-2\rho}f_k(s)ds. 
	$$ 

Given this stochastic stability result, one may consider stabilization using multiplicative noise on systems that are apriori modelled using sector bounded uncertainties, or systems that can be approximated as such.  The effect of SAR is quite dramatic as ca be seen in Fig. 1 and Fig. 2, from a recent publication of the authors \cite{syent}, where multiplicative noise is shown to stabilize a noise-free chaotic behavior of a non linear system \cite{Kwok} with sector bounded nonlinearities.

\section{Application to More General Systems}

We note that the model of system (\ref{e1bis}) and (\ref{sector1}) is relevant also in cases where $f (y(t))$ in the model is not a priori sector-bounded. In such
cases, one may invoke the universal approximation theorem \cite{Cyb} to systems where a single hidden layer, with, e.g., a tanh activation function and a linear output layer, provides an approximation with arbitrarily small error for an arbitrarily wide hidden layer. In such
cases, the model of system (4) readily becomes relevant, as the approximate function is now sector-bounded. However, one should be very careful, as systems exhibiting chaotic behaviors may require a very high degree of approximation to maintain their chaotic nature.

One such example is the Morris-Lecar model of a neuron with three Ion channels and an external applied current $I$. The details for this model, were taken from \cite{Moris}.
The model comprises of a second-order nonlinear differential equation

$$
dV/dt = I - g_L (V - V_L) - g_Ca  M_ss (V - V_{Ca}) - g_K  N (V - V_K))/C
$$
and
$$
dN/dt = (N_{ss} - N) / \tau_N
$$
where
$$
M_{ss} = \frac{1}{2}(1 + tanh((V - V_1) / V_2))
$$
$$
N_{ss} = \frac{1}{2}(1 + tanh((V - V_3) / V_4))
$$

$$
\tau_N = 1 / (\phi cosh((V - V_3) / (2 V_4)))
$$
in which $V$ is the voltage across the membrane, $N$ is the recovery variable representing the instantaneous probability that the $K^{+}$ ion channel is in its conducting state and where $C$ is the capacitance of the cell membrane. The values of the constant parameters in this model, as well as the input current $I$, determine the interplay between the voltage $V$ and the recovery variable $N$ and to lead to either convergence to a steady state or oscillations.     

Those different type of behaviors, can be observed in the responses of the three ion channels, $L^{+}$, $C_a^{2+}$ and $K^{+}$ respectively.
$$
L = - g_L  (V - V_L)
$$
$$
Ca =  - g_{Ca}  M_{ss}  (V - V_{Ca})
$$
$$
K = - g_K  N  (V - V_K)
$$
Setting the tuning parameters to $C = 5, V_1 = -1.2, V_3 = 12, V_4 = 17.4, \phi = 1/15, V_L = -60, V_{Ca} = 120, 
V_K = -80, g_Ca = 4, g_K = 8, g_L = 2$ and the initial conditions to $V = -52.14$ and $N = 0.02$ we obtain the oscillatory behavior of Fig. 3. We next try to stabilize the oscillations, by injecting a multiplicative noise component in $I$. 


We define the following functions that are, apparently, non sector bounded, but in order to allow SAR analysis we approximate them, using a shallow neural network that utilizes  sector bounded activation functions, 
$$
f_i(x) =  W_2^i ( tanh( W_1^i  x +  b_1^i ))+ b_2^i , i=1,2,3
$$
The weight matrices and bias terms are found by backpropagation training (see e.g. \cite{Haykin}). In this case, we approximated each of the $f_i$ using a network with a single hidden layer and $10$ neurons. Therefore, the functions vector $f$ of (\ref{e1bis}) is of order $30$ and accordingly $F$ is of dimension $2 \times 30$. To fit the resulting model to (\ref{e1bis}) we add $28$ fictititous states which are zero. To this end, we denote for $\kappa >0$, 
$$
\bar A = \left[ \begin{array}{cc} A & 0 \\ 0 & -\kappa I \end{array} \right],  F = \left[ \begin{array}{c} F \\ 0 \end{array} \right]
$$
Also $S$ and $\Delta$ are derived from $W_1^i, i=1,2,3$. Such an approximation allows using the above stochastic stability analysis, to search and verify the closed-loop stability of the Morris-Lecar neuron, under the SAR controller. Indeed, in Fig. 5, a sweep of SAR level $\sigma$ values to check the $\mathcal N<0$ condition of (\ref{e18}) using YALMIP (\cite{yalmip}). The simulation results of Fig. 4 with $\sigma= 0$ and $\sigma = 0.85$ agree with those LMI based stochastic stability analysis results. 

One should note that the result of Fig. 5  is quite noisy when SAR is applied, and not as strict nulling of output as in the case of Fig. 1 and Fig. 2. and to ease viewing the stabilisation effect of SAR, filtered version of $V$ and $N$ have been displayed as well as the original state vector components of the Moris-Lecar model. The noisy result in the non-filtered states, motivates the future work of adding performance maeasures rather just than stability, such the $H_2$ norm of selected objective, e.g. covariance of $C_1 E\{ xx^T\}C_1^T$.

\section{Concluding Remarks}
Results regarding stochastic stability of SAR controllers applying state-multiplicative noise have been reviewed. 
The universal approximation theorem, considerably widens the scope of those results, and, indeed, we have shown application to the Morris-Lecar model which is not apriori described in terms of sector-bounded functions. Using an approximation with a neural network with a single hidden layer and $10$ $tanh$ activated neurons, it was shown that the LMI based stability analysis, correctly predicts the SAR gain (i.e. the gain on state-multiplicative noise ) range for stability. It should be noted that further applications can be obtained following \cite{syent2}, where mixed controllers applying a combination of stochastic and deterministic gain.


\newpage

\begin{figure}
	\centering
	{\epsfig{file = 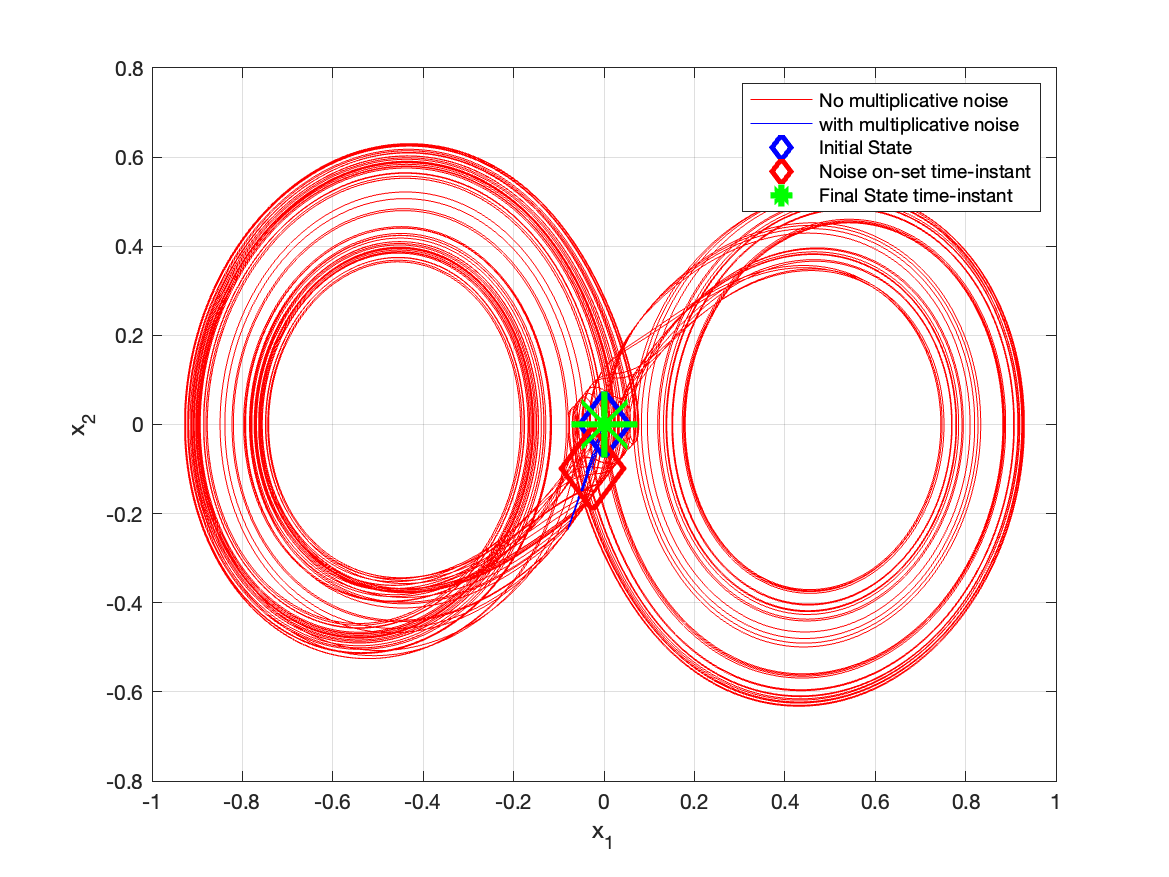, width = 12cm}}
	\caption{$x_1$-vs. $x_2$ Sector Bounded Nonlinearities}
\end{figure}

\begin{figure}[h]
	\centering
	\includegraphics[width=8cm, height=8cm]{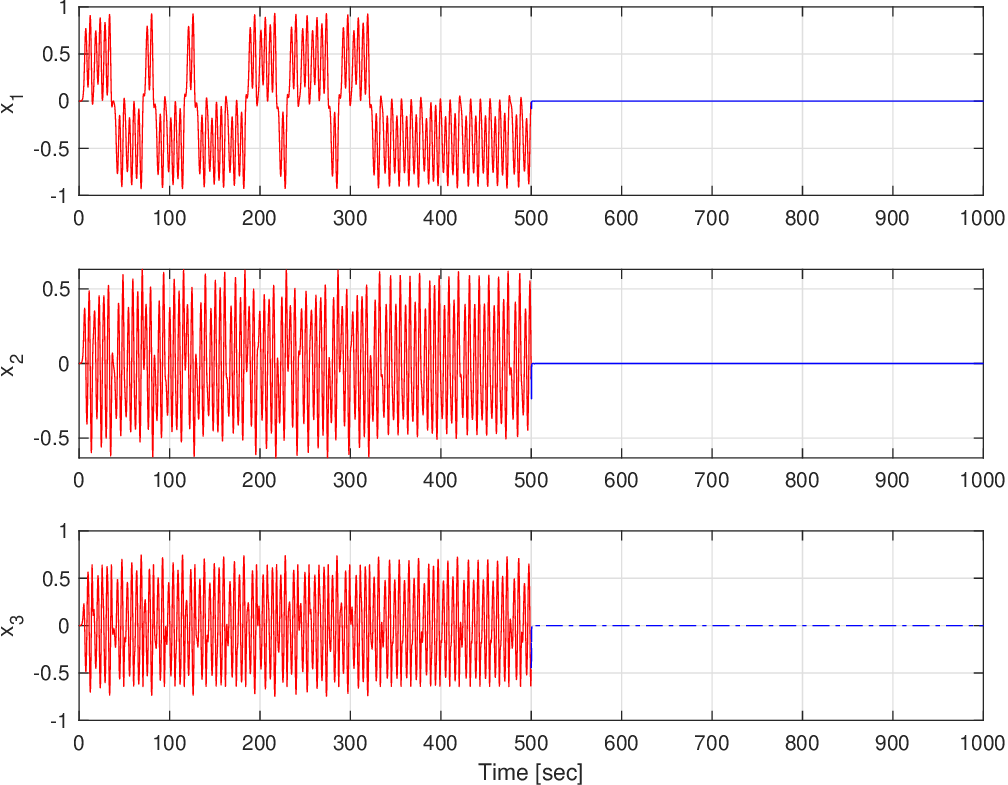}
	\caption{$x_1(t)$ and $x_2(t)$ Sector Bounded Nonlinearities}
\end{figure}

\begin{figure}
	\centering
	{\epsfig{file = 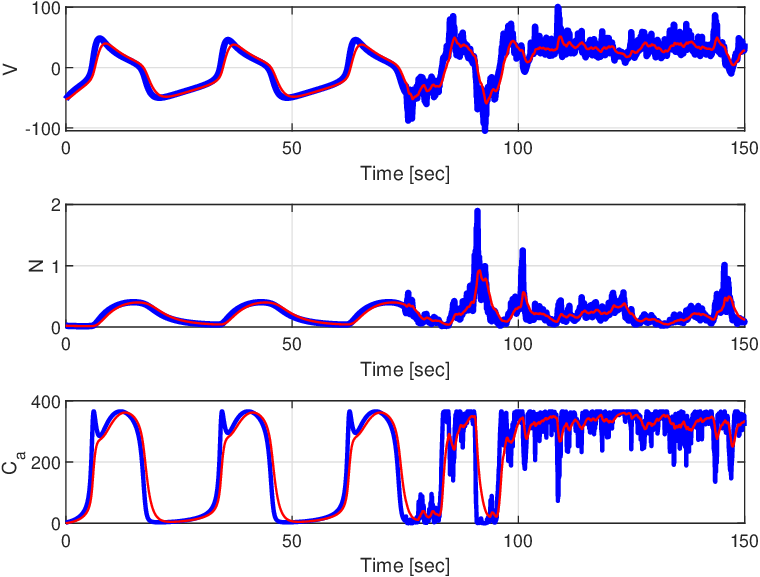, width = 12cm}}
	\caption{$x_1(t)$ and $x_2(t)$ Moris-Lecar Model - with multiplicative noise}
\end{figure}

\begin{figure}
	\centering
	{\epsfig{file = 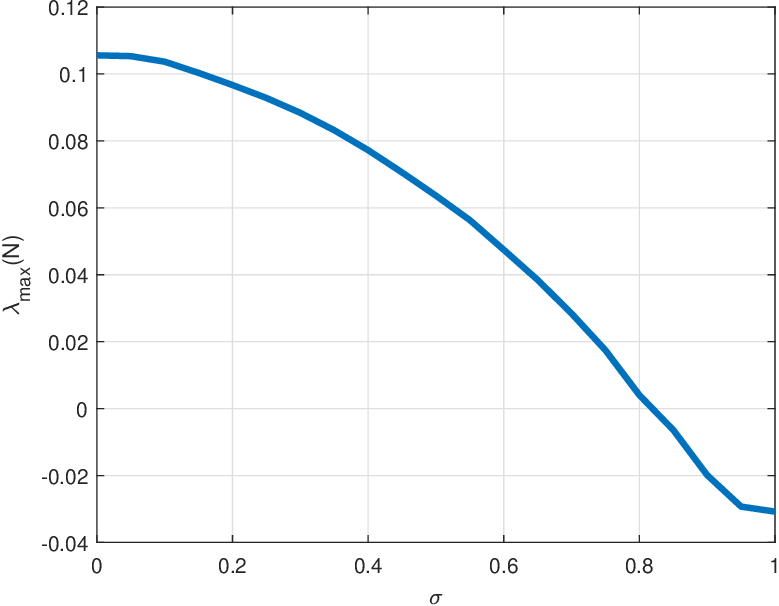, width = 12cm}}
	\caption{Sweep of $\sigma$ for $\mathcal N<0$ condition}
\end{figure}

\enddocument